\documentclass[11pt]{article}
\usepackage{amsmath}
\usepackage{amssymb}
\usepackage{amsfonts}
\usepackage{amsthm}
 \usepackage{amsthm}
 \newtheorem{thm}{Theorem}
 \theoremstyle{plain}
 \newtheorem{prop}[thm]{Proposition}
 \theoremstyle{plain}

\newcommand{\R}{\mathbb{R}}
\begin{document}
%
%
%
\title{Fundamental Flaws in Feller's\\
Classical Derivation of Benford's Law}
\author{Arno Berger\footnote{Arno Berger is Associate Professor with the
Department of Mathematical and Statistical Sciences, University of Alberta, Edmonton, Alberta T6G 2G1, Canada (email:
{\tt aberger@math.ualberta.ca}); he was supported by an NSERC Discovery Grant.}\\[1mm]
{\small Mathematical and Statistical Sciences, University of Alberta}\\[4mm]
and\\[4mm]
Theodore P.\ Hill\footnote{Theodore P.\ Hill is Professor Emeritus with the
School of Mathematics, Georgia Institute of Technology, Atlanta GA 30332-0160, USA (email: {\tt hill@math.gatech.edu}).
}\\[1mm]
{\small School of Mathematics, Georgia Institute of Technology}\\[10mm]}
\date{{\small 14 May 2010}}
\maketitle
\begin{abstract}
\noindent
Feller's classic text {\em An Introduction to Probability Theory and its 
Applications\/} contains a derivation of the well known significant-digit law 
called Benford's law. More specifically, Feller gives a sufficient condition 
(``large spread'') for a random variable $X$ to be approximately Benford distributed, 
that is, for $\log_{10}X$ to be approximately uniformly distributed modulo one. 
This note shows that the large-spread derivation, which continues to be widely 
cited and used, contains serious basic errors. Concrete examples and a new 
inequality clearly demonstrate that large spread (or large spread on a logarithmic 
scale) does not imply that a random variable is approximately Benford distributed, 
for any reasonable definition of ``spread'' or measure of dispersion.
\end{abstract}

\vspace{1ex} \noindent {\bf Key words:} Significant digit, 
Benford distribution, Kolmogorov--Smirnoff distance, uniform distribution modulo one.
 
\clearpage

\noindent
In probability and statistics, a correct general explanation of a 
principle is often as valuable as a detailed formal argument. In his December 
2009 column in the IMS Bulletin, UC Berkeley statistics professor T.\ Speed 
extols the virtues of derivations in statistics (\cite{S}):

\begin{quote}
{\em I think in statistics we need derivations, not proofs. That is, lines of 
reasoning from some assumptions to a formula, or a procedure, which may 
or may not have certain properties in a given context, but which, all 
going well, might provide some insight.}
\end{quote}

\noindent
For illustration, Speed quotes two examples of the convolution property 
for the Gamma and Cauchy distributions from the classic 1966 text {\em An 
Introduction to Probability Theory and Its Applications\/} by W. Feller (\cite{Fel}).

On page 63, Feller also gave a brief derivation, in Speed's sense, of the 
well known logarithmic distribution of {\em significant digits\/} 
called Benford's law (\cite{Ben, Few, H1, H2, N, R}).
Recall that if a random variable is Benford (i.e.\ has a Benford distribution) 
then its first significant digit is ``1'' with probability
$\log_{10}2 \approx 0.3010$; similar expressions hold for the general joint Benford 
distributions of all the significant digits (\cite{H1}). For the purposes of this note, 
a simple and very useful characterization of a Benford distribution is 
\bigskip

\noindent
{\bf (1)} A positive random variable $X$ is Benford if and only if
$\log_{10}X$ is uniformly distributed (mod 1).

\bigskip

\noindent
Since Feller has inspired so many who teach probability and statistics today, 
and since many undergraduate courses now include a brief introduction to 
Benford's law, it is not surprising that Feller's derivation is still in 
frequent use to provide some insight about Benford's law. For example, 
a class project report for a 2009 upper-division course in statistics at UC
Berkeley (\cite[p.3]{AP1}) said

\begin{quote}
{\em \dots like the birthday paradox, an explanation [of Benford's law] occurs 
quickly to those with appropriate mathematical background \dots To a mathematical 
statistician, Feller's paragraph says all there is to say \dots Feller's derivation 
has been common knowledge in the academic community throughout the last 40 years.}
\end{quote}

\noindent
The online database \cite{BH} lists about twenty published references since 2000 
alone to Feller's argument (e.g.\ \cite{AP1, Few}) the crux of which is 
Feller's claim (trivially edited) that 

\bigskip

\noindent
{\bf (2)} If the spread of a random variable $X$ is very large, then $\log_{10} X$ (mod 1) 
will be approximately uniformly distributed on $[0, 1)$.

\bigskip

\noindent
The implication of (1) and (2) is that all random variables with large spread 
will be approximately Benford distributed. That sounds quite plausible, but as 
C.S.\ Pierce observed (\cite[p.174]{Gar}), ``in no other branch 
of mathematics is it so easy for experts to blunder as in probability theory''. 
Indeed, even Feller blundered on Benford's law, and took many experts with him. Claim 
(2) is simply false under any reasonable definition of spread or measure 
of dispersion, including {\em range}, {\em interquartile range\/} (or distance between the
$(1-\alpha)$- and the $\alpha$-quantile), {\em standard deviation}, or {\em mean difference\/} 
(Gini coefficient), no matter how smooth or level a density the random variable $X$ may have. 
To see this, one does not have to look far. Concretely, no positive uniformly distributed 
random variable even comes close to being Benford, regardless of how large (or small) its 
spread is. This statement can be quantified explicitly via the following new inequality; 
for its formulation, recall that the Kolmogorov-Smirnoff distance $d_{\rm KS}(X,Y)$ between 
two random variables $X$ 
and $Y$ with cumulative distribution functions $F$ and $G$, respectively, is $d_{\rm KS}(X,Y)=
\sup \{ |F(x) - G(x)| : x\in \R\}$.

\begin{prop}[\cite{Ber}]
For every positive uniformly distributed random variable $X$, 
$$
d_{\rm KS} \bigl(
\log_{10} X \, \mbox{\rm (mod}\, 1), \, U(0,1)
\bigr) \ge
\frac{-9 + \ln 10 + 9 \ln 9 - 9 \ln \ln 10}{18 \ln 10} = 0.1334 \ldots \, ,
$$
and this bound is sharp.
\end{prop}

There is nothing special about the usage of the Kolmogorov-Smirnoff distance or decimal 
base in this regard; similar universal bounds hold for the Wasserstein distance, for 
example, and other bases. Another way to see that (2) is false, in the discrete and 
significant-digit setting, is to observe that no matter how large $n$ is, an integer-valued 
random variable uniformly distributed on the first $2\cdot 10^n$ positive integers will have more than 50\% 
of its values beginning with a ``1'', as opposed to the Benford probability of about 30\%.

How could Feller's error have persisted in the academic community, among students and 
experts alike, for over 40 years? Part of the reason, as one colleague put it, is simply 
that Feller, after all, is Feller, and Feller's word on probability has just been
taken as gospel. Another reason for the long-lived propagation of the error has apparently 
been the confusion of (2) with the similar claim

\bigskip

\noindent
{\bf (3)} If the spread of a random variable $X$ is very large, then $X$ (mod 1) will be 
approximately uniformly distributed on $[0, 1)$.

\bigskip
\noindent
For example, \cite[p.3]{AP1} cites Feller's claim (2), but \cite[p.8]{AP1}
cites Feller's claim as (3). A third possible explanation for the persistence of the error is the 
common assumption that (3) implies (2). For example, \cite[p.1]{GD} state: 

\begin{quote}
{\em An elementary new explanation has recently been published, based on the fact that any
$X$ whose distribution is ``smooth'' and ``scattered'' enough is Benford. The 
scattering and smoothness of usual data ensures that $\log (X)$ is itself smooth 
and scattered, which in turn implies the Benford characteristic of $X$.}
\end{quote}

\noindent
Now (3) is also intuitive and plausible, but unlike (2), it is often accurate
if the distribution is fairly uniform. And if the distribution is not fairly uniform, 
then without further information, no interesting conclusions at all can be made about 
the significant digits --- most of the values could for instance start with a ``7''. Since $X$ has very, 
very large spread if and only if $\log X$ has very large spread, on the surface (2) 
and (3) appear to be equivalent. After all, what difference can one tiny extra ``very'' 
mean? On the other hand, as Proposition 1 clearly implies, they are not the same, 
and (2) is false.

Although (3) is perhaps more accurate than (2), unfortunately it does not explain Benford's 
law at all, since the criterion in (1) says that $X$ is Benford if and only if the 
logarithm of $X$ --- and not $X$ itself --- is uniformly distributed (mod 1). 
Some authors partially explain the ubiquity of Benford distributions based on 
an assumption of a ``large spread on a logarithmic scale'' (e.g.\ \cite{AP1, AP2, Few, W}).
Others (e.g.\ \cite[p.17]{AP2}) claim that ``what Feller obviously {\em meant}'' [italics in original] 
by spread was log spread, i.e. that when Feller wrote (2) he really meant to say that

\bigskip
 
\noindent
{\bf (3')} If $\log_{10} X$ has very large spread, then $\log_{10} X$ (mod 1) will be 
approximately uniformly distributed on $[0,1)$,

\bigskip

\noindent
which is but an unnecessarily convoluted version of (3).
They then apply (3) or (3') to conclude that if $\log_{10}X$ has large spread, 
then $X$ is approximately Benford. This avoids Feller's error (2), but still leaves 
open the question of why it is reasonable to assume that the {\em logarithm\/} of the spread, 
as opposed to the spread itself or, say, the $\log \log$ spread should be large. 
As seen above, those assumptions contain a subtle difference, and lead to very 
different conclusions about the distributions of the significant digits. 
Using the same logic, for instance, an assumption of large spread on the log log
scale would imply that $\log X$ is Benford, whereas none of the usual Benford
random variables such as $X_k$ with densities $1/(x \ln 10)$ on $(10^k, 10^{k+1})$
are also Benford on the log scale. Moreover, 
via (1) and (3), assuming large spread on a logarithmic scale is equivalent to 
assuming an approximate Benford distribution. Quite likely, Feller realized this, 
and in (2) specifically did {\em not\/} hypothesize that the $\log$ of the range was large. 

A related and apparently widespread misconception is that claim (2),
notwithstanding its incorrectness, or claim (3) implies that a larger spread
or log spread automatically means closer conformance to Benford's law. 
For example, \cite{W} concludes that ``datasets with large logarithmic spread 
will naturally follow the law, while datasets with small spread will not'', 
and the Conclusion of the study \cite[p.12]{AP2} states

\begin{quote}
{\em On a small stage (18 data-sets) we have checked a theoretical prediction. 
Not just the literal assertion of Benford's law - that in a data-set with large 
spread on a logarithmic scale, the relative frequencies of leading digits will 
approximately follow the Benford distribution - but the rather more specific 
prediction that distance from Benford should decrease as that spread increases. 
In one sense it's not surprising this works out.}
\end{quote}

\noindent
But distance from the Benford distribution does not generally decrease 
as the spread increases, regardless of whether the spread is measured 
on the original scale or on the logarithmic scale. A simple way to see 
this is as follows: Let $Y$ be a random variable uniformly distributed 
on $(0,1)$, and let $X=10^Y$ and $Z=10^{3Y/2}$. Then by (1), $X$ is exactly 
Benford, since $\log_{10}X=Y$, and $Z$ is not close to Benford since $3Y/2$ 
(mod 1) is not close to uniform on $(0,1)$. Yet for any reasonable definition 
of spread, including all those mentioned earlier, the spread of $Z$ is 
larger than the spread of $X$, and the spread of $\log_{10}Z = 3Y/2$ is larger 
than the spread of $\log_{10}X = Y$. Another way to see that the distance 
from the Benford distribution does not decrease as the spread increases 
is contained in the proof of Proposition 1: For $X_T$ a random variable 
uniformly distributed on $(0,T)$, it is shown that the Kolmogorov-Smirnoff 
distance between $\log_{10}X_T$ and $U(0,1)$ is a continuous $1$-periodic function of $\log_{10}T$.
Moreover, when employing a logarithmic scale it is important to keep in mind that
what is considered large generally depends on
the base of the logarithm. For example, as noted earlier, if $Y$ is uniformly distributed
on $(0,1)$ then $X=10^Y$ is exactly Benford base $10$, yet it is not Benford base $2$
even though its spread on the $\log_2$-scale is $\log_2 10\approx 3.3219$ times
as large.

It is interesting to note that when Feller credited Pinkham in his derivation 
in 1966, it was not widely known that Pinkham's argument (\cite{P})
for the scale-invariant characterization of Benford's law also 
contains an irreparable and fundamental flaw. Raimi (\cite[sec.7]{R})
explains Pinkham's error in detail, and credits Knuth (\cite{K})
for the discovery that the error was in Pinkham's unwritten implicit assumption 
that there exists a scale-invariant probability distribution on the positive real 
numbers --- when clearly there does not, since the largest median of every positive 
random variable changes under changes of scale. The first correct proof that 
the Benford distribution is the unique scale-invariant probability distribution
on the significant digits (and the unique continuous 
base-invariant distribution) is in \cite{H2}.

In conclusion, classroom experiments based on Feller's derivation or on an assumption of 
large range on a logarithmic scale (e.g.\ \cite{AP1, AP2, Few, W}) 
should be used with caution. As an alternative or supplement, teachers might also ask students 
to compare the significant digits in the first 20-30 articles in tomorrow's {\em New 
York Times\/} against Benford's law, thereby testing real-life data against the 
explanation given in the main theorem in \cite{H2} which, without any assumptions
on magnitude of spread, shows that mixing data from different 
distributions in an unbiased manner leads to a Benford distribution.

Although some experts may still feel that ``like the birthday paradox, there is a 
simple and standard explanation'' for Benford's law (\cite[p.6]{AP2}) and 
that this explanation occurs quickly to those with appropriate 
mathematical background, there does not appear to be a simple derivation of Benford's law 
that both offers a ``correct explanation'' (\cite[p.7]{AP2}) and satisfies Speed's 
goal to provide insight. In that sense, although Benford's law now rests on 
solid mathematical ground, most experts seem to agree with \cite{Few} that its ubiquity in 
real-life data remains mysterious.

\medskip

\noindent
{\bf Acknowledgement.} The authors are grateful to Rachel Fewster, Kent Morrison, 
and Stan Wagon for several helpful communications.


\begin{thebibliography}{123}


{\small

\bibitem[AP1]{AP1}
Aldous, D., and Phan, T.\ (2009), When Can One Test an Explanation? 
Compare and Contrast Benford's Law and the Fuzzy CLT, Class project report dated May 11, 2009, 
Statistics Department, UC Berkeley; accessed on May 14, 2010 at \cite{BH}.

\bibitem[AP2]{AP2}
Aldous, D., and Phan, T.\ (2010), When Can One Test an Explanation? 
Compare and Contrast  Benford's Law and the Fuzzy CLT, Preprint dated Jan.\ 3, 2010, 
Statistics Department, UC Berkeley; accessed on May 14, 2010 at \cite{BH}.

\bibitem[Ben]{Ben}
Benford, F.\ (1938), The law of anomalous numbers, {\em Proceedings of the 
American Philosophical Society\/} {\bf 78}, 551--572. 

\bibitem[Ber]{Ber}
Berger, A.\ (2010), Large spread does not imply Benford's law, Preprint;
acessed on May 14, 2010 at {\tt http://www.math.ualberta.ca/$\sim$aberger/Publications.html}.

\bibitem[BH]{BH}
Berger, A., and Hill, T.P.\ (2009), {\em Benford Online Bibliography}, 
{\tt http://www.benfordonline.net}; accessed May 14, 2010. 

\bibitem[Fel]{Fel}
Feller, W.\ (1966), {\em An Introduction to Probability Theory and Its Applications\/} vol 2, 
2nd ed., J.\ Wiley, New York.

\bibitem[Few]{Few}
Fewster, R.\ (2009), A simple Explanation of Benford's Law, {\em American Statistician} {\bf 63}(1), 20--25. 

\bibitem[Ga]{Gar}
Gardner, M.\ (1959), Mathematical Games: Problems involving questions of probability and
ambiguity, {\em Scientific American} {\bf 201}, 174--182.

\bibitem[GD]{GD}
Gauvrit, N., and Delahaye, J.P.\ (2009), Loi de Benford g\'en\'erale, 
{\em Math\'ematiques et sciences humaines\/} {\bf 186}, 5--15;
accessed May 14, 2010 at {\tt http://msh.revues.org/document11034.html}. 

\bibitem[H1]{H1}
Hill, T.P.\ (1995), Base-Invariance Implies Benford's Law, {\em Proceedings of the 
American Mathematical Society\/} {\bf 123}(3), 887--895.

\bibitem[H2]{H2}
Hill, T.P.\ (1995), A Statistical Derivation of the Significant-Digit Law, {\em Statistical Science\/} 
{\bf 10}(4), 354--363.

\bibitem[K]{K}
Knuth, D.\ (1997), {\em The Art of Computer Programming}, pp 253-264, vol. 2, 3rd ed, 
Addison-Wesley, Reading, MA.

\bibitem[N]{N}
Newcomb, S.\ (1881), Note on the frequency of use of the different digits in natural numbers,
{\em American Journal of Mathematics\/} {\bf 4}(1), 39--40.

\bibitem[P]{P}
Pinkham, R.\ (1961), On the Distribution of First Significant Digits,
{\em Annals of Mathematical Statistics\/} {\bf 32}(4), 1223--1230.

\bibitem[R]{R}
Raimi, R.\ (1976), The First Digit Problem, {\em American Mathematical Monthly} {\bf 83}(7), 521--538.

\bibitem[S]{S}
Speed, T.\ (2009), You want proof?, {\em Bulletin of the Institute of Mathematical Statistics\/} 
{\bf 38}, p 11.

\bibitem[W]{W}
Wagon, S.\ (2010), Benford's Law and Data Spread; accessed May 14, 2010 
at Wolfram Online Demonstrations Projects 
{\tt http://demonstrations.wolfram.com/BenfordsLawAndDataSpread}.

}

\end{thebibliography}
\end{document}